\documentclass[12pt]{article}
\usepackage{amsmath}
\usepackage{amsfonts,cite}
\usepackage{mathrsfs}
\usepackage{amssymb}
\usepackage{amsthm}
\usepackage{graphicx}
\usepackage{epsfig}
\usepackage{epstopdf}

\renewcommand\baselinestretch{1.2}
\allowdisplaybreaks

\title{\textbf{Rigidity of spacelike translating
 solitons in pseudo-Euclidean space}
\thanks{\baselineskip 9pt Research partially
supported by NSFC.}}
\author{{Ruiwei Xu \quad\quad Tao Liu}}

\date{}
\begin{document}
\maketitle
\renewcommand{\baselinestretch}{1}
\renewcommand{\arraystretch}{1.2}

\date{}
\maketitle

\date{}
\catcode`@=11 \@addtoreset{equation}{section} \catcode`@=12
\maketitle{}
\makeatletter
\renewcommand{\theequation}{\thesection.\arabic{equation}}
\@addtoreset{equation}{section} \makeatother

{\bf Abstract:} In this paper, we investigate the parametric version and non-parametric version of rigidity theorem of spacelike translating solitons
in pseudo-Euclidean space $\mathbb{R}^{m+n}_{n}$.  Firstly, we classify $m$-dimensional complete spacelike translating solitons in $\mathbb{R}^{m+n}_{n}$ by affine technique and classical gradient estimates, and prove the only complete spacelike translating solitons in $\mathbb{R}^{m+n}_{n}$ are the spacelike $m$-planes. This result provides another proof of a nonexistence theorem for complete spacelike translating solitons in \cite{C-Q}, and a simple proof of rigidity theorem in \cite{X-H}. Secondly, we generalize the rigidity theorem of entire spacelike Lagrangian translating solitons in \cite{X-Z} to spacelike translating solitons with general codimensions. As a directly application of theorem, we obtain two interesting  corollaries in terms of Gauss image.


\vskip 6pt  {\bf 2000 AMS Classification:}
53C44; 53C40; 53C43.

\section{Introduction}

The mean curvature flow (MCF) in Euclidean space (pseudo-Euclidean space resp.) is a one-parameter family of immersions
$$
X_{t}=X(\cdot,t):M^{m}\rightarrow \mathbb{R}^{m+n}(\mathbb{R}^{m+n}_{n}\,\, resp.)
$$
with the corresponding image $M_{t}=X_{t}(M)$
such that
\begin{equation}\label{1.1}
\left\{
\begin{array}{lll}
\frac{d}{dt}X(x,t)=H(x,t), ~~~ x\in M,\\
X(x,0)=X(x),\\
\end{array}
\right.
\end{equation}
is satisfied, here $H(x,t)$ is the mean curvature vector of $M_{t}$ at $X(x,t)$ in $\mathbb{R}^{m+n}$ ($\mathbb{R}^{m+n}_{n}$ resp.). The MCF in higher codimension has been studied extensively in the last few years (cf.\cite{Neves,C-J-Y-2,S-K-2,X-Y-L-1,X-Y-L-2}). Translating soliton and self-shrinker play an important role in analysis of singularities in MCF. A submanifold $M^{m}$ is said to be a translating soliton in $\mathbb{R}^{m+n}$ (spacelike translating soliton in $\mathbb{R}^{m+n}_{n}$ resp.) if the mean curvature vector $H$ satisfies
\begin{equation}\label{1.2}
H=T^{\bot},
\end{equation}
where $T\in\mathbb{R}^{m+n}$ is a non-zero constant vector, which is called a translating vector, and $T^{\bot}$ denotes the orthogonal projection of $T$ onto the normal bundle of $M$. $M^{m}$ is said to be a self-shrinker in $\mathbb{R}^{m+n}$(spacelike self-shrinker in $\mathbb{R}^{m+n}_{n}$ resp.) if it satisfies a quasi-linear elliptic system
\begin{equation}\label{1.3}
H=-X^{\bot},
\end{equation}
where $X^{\bot}$ is the normal part of $X$, which is an important class of solutions to (\ref{1.1}).

There is a plenty of works on the classification and uniqueness problem for translating soliton and self-shrinker in Euclidean space
(cf. \cite{B-S,K-K,T-H,A-N,H-G,S-K-3,C-L,D-X-Y,W-L,C-P}).
On the other hand, there are many works on the rigidity problem for complete spacelike submanifolds in pseudo-Euclidean space. Calabi \cite{E-C} proposed the Bernstein problem for spacelike extremal hypersurfaces in Minkowski space $\mathbb{R}^{m+1}_{1}$ and proved such hypersurfaces have to be hyperplanes when $m\leq 4$. In \cite{C-S-Y} Cheng-Yau solved the problem for all dimensions. Later, Jost-Xin \cite{J-X} generalized the results to higher codimensions.

\vskip 0.1in  \textbf{Theorem 1 \cite{J-X}:}
\emph{Let $M^m$ be a spacelike extremal submanifold in $\mathbb{R}^{m+n}_{n}$. If $M^m$ is closed with respect to the Euclidean topology, then $M$ has to be a linear subspace.}

Translating soliton and self-shrinker can be regarded as generalizations of extremal submanifold. Thus it is natural to study the corresponding rigidity problem for spacelike translating soliton and self-shrinker.  Here we only investigate the rigidity of spacelike translating soliton in $\mathbb R^{m+n}_n$.

For complete spacelike Lagrangian translating solitons in $\mathbb R^{2n}_n$, Xu-Huang \cite{X-H} proved the following rigidity theorem.

\vskip 0.1in  \textbf{Theorem 2\cite{X-H}:}\emph{
Let $f(x_1,...,x_n)$ be a strictly convex $C^{\infty}$-function defined on a convex domain $\Omega\subset\mathbb{R}^n$. If the graph $M_{\nabla f}=\{(x,\nabla f(x))\}$ in $\mathbb R^{2n}_n$ is a complete spacelike translating soliton, then $f(x)$ is a quadratic polynomial and $M_{\nabla f}$ is an affine $n$-plane.}

\vskip0.1in

From Theorem 2 above, it is easy to see that the corresponding translating vector must be \emph{\textbf{spacelike}}. In \cite{C-Q}, Chen-Qiu proved a nonexistence theorem for complete spacelike translating solitons in $\mathbb{R}^{m+n}_{n}$ by establishing a very powerful generalized Omori-Yau maximum principle.
They proved that \emph{there exists no complete $m$-dimensional spacelike translating soliton (with a \textbf{timelike} translating vector)}. In fact, such nonexistence conclusion still holds for the lightlike translating vector case. Motivating by papers\cite{C-Q,X-H}, we classify $m$-dimensional complete spacelike translating solitons in $\mathbb{R}^{m+n}_{n}$ by affine technique and classical gradient estimates, and obtain the following Bernsterin theorem.

\vskip 0.1in  \textbf{Theorem 3:}\emph{
Let $M^m$ be a complete spacelike translating soliton in $\mathbb{R}^{m+n}_{n}$, then it is an affine $m$-plane.}

\vskip 0.1in
A more precise statement of the assertion in Theorem 3 says that there exists an $m$-dimensional complete spacelike translating soliton in $\mathbb{R}^{m+n}_{n}$ only if the translating vector is \textbf{spacelike}. It provides another proof of the nonexistence theorem of translating soliton in \cite{C-Q}, and a simple proof of the Bernstein theorem for translating soliton in \cite{X-H}. Here we mention that it is valid if one shall use the generalized Omori-Yau maximum principle in \cite{C-Q} to prove Theorem 3 above.

For non-parametric Lagrangian solitons of mean curvature flow, there are several interesting rigidity theorems. As to entire self-shrinker for mean curvature flow in $\mathbb{R}^{2n}_{n}$ with the indefinite metric $\sum_{i=1}^{n}dx_idy_i$, Huang-Wang \cite{H-W} and Chau-Chen-Yuan \cite{C-C-Y} used different methods
to prove the rigidity of entire self-shrinker under \emph{a decay condition on the induced metric $(D^2f)$}.
Later, using an integral method, Ding and Xin \cite{D-X} removed the additional decay condition and proved
\emph{any entire smooth convex self-shrinking solution for mean curvature flow in $\mathbb{R}^{2n}_{n}$ is
a quadratic polynomial.} In \cite{H-X}, Huang and Xu investigated the rigidity problem of entire spacelike translating soliton graph $(x,\nabla f)$ in $\mathbb{R}^{2n}_{n}$ with some symmetry conditions. Motivating by\cite{C-C-Y,D-X,H-W,H-X}, Xu-Zhu\cite{X-Z} proved a rigidity  theorem of entire convex translating solutions for mean curvature flow in $\mathbb{R}^{2n}_{n}$ under a decay condition on the induced metric $(D^2f)$ and provided a class of nontrivial entire spaclike Lagrangian translating solitons.

\vskip 0.1in  \textbf{Theorem 4\cite{X-Z}:}
Let $f(x)$ be an entire smooth strictly convex function on $\mathbb R^n$ $(n\geq2)$ and its graph $M_{\nabla f}=\{(x,\nabla f(x))\}$ be a translating soliton in $\mathbb R^{2n}_n$. If there exists a number $\epsilon>0$ such that the induced metric $(D^{2}f)$ satisfies
\begin{equation}
(D^2f)>\frac{\epsilon}{|x|^{2}}I,\;\,\,as\,\,\, |x|\rightarrow\infty , \end{equation}
then  $f(x)$ must be  a quadratic polynomial and $M_{\nabla f}$ is an affine $n$-plane.

\vskip 0.1in

Here we shall use the idea in \cite{X-Z} to generalize Theorem 4 to spacelike graphic translating solitons in pseudo-Euclidean space $\mathbb{R}^{m+n}_n$.

\vskip 0.1in  \textbf{Theorem 5:}
Let $u^\alpha (1\leq\alpha\leq n)$ be smooth functions defined everywhere in $\mathbb{R}^m$ and their graph $M=(x,u^1(x),u^2(x),\cdots,u^n(x))$ be a spacelike translator in $\mathbb{R}^{m+n}_n$.  If there exists a number $\epsilon>0$ such that  the induced metric $(g_{ij})$ satisfies
\begin{equation}
(g_{ij})>\frac{\epsilon}{|x|}I,\;\;\;as\,\,\,\;\; |x|\rightarrow\infty , \end{equation}
then $u^{1}(x),\cdots,u^{n}(x)$ are linear functions on $\mathbb{R}^m$, and  $M$ is an affine $m$-plane in $\mathbb{R}^{m+n}_n$.

\vskip 0.1in

\textbf{Remark}: It is necessary that there is a restriction on the induced metric $(g_{ij})$ for the rigidity of translating solitons. If not, there exist nontrivial entire smooth spacelike translating solitons, which are not planes.
For example, submanifold
\begin{equation}\label{e1.5}
\left(x,y,\ln(1+\exp\{2x\})-x, \mu y\right),\quad |\mu|< 1,\quad (x,y)\in\mathbb R^2
\end{equation}
is an entire spacelike translator with the translating vector $(0,1,1,\mu)$ in $\mathbb{R}_{2}^{4}$, which graphic functions satisfy the PDE (\ref{2.12}).

By studying the distribution of the Gauss map, they obtained Bernstein theorems of translating solitons in Euclidean space (see \cite{B-S},\cite{K-K} and\cite{X-Y-L-5}). Notice that the Gauss image of spacelike graphic submanifold $M^m$ in $\mathbb{R}^{m+n}_n$ is bounded if and only if the induced metric $g=\det(g_{ij})$  is bounded (see \cite{X-Y-L-4}). Therefore, it is easy to see that example (1.6) above and example (1.7) in \cite{X-Z} have boundless Gauss images. As a directly application of theorem 5, we have

\vskip 0.1in  \textbf{Corollary 1:}
Let $M^m$ be an entire spacelike graphic translator in $\mathbb{R}^{m+n}_n$ as defined in Theorem 5. If the Gauss image of $M$ is bounded,
then $M$ is an affine $m$-plane.

\vskip 0.1in By relaxing the bound of the Gauss image to controlled growth, we also get a more general corollary from theorem 5 as Prof. Dong generalized a rigidity theorem for spacelike graph with parallel mean curvature in\cite{D}.

\vskip 0.1in  \textbf{Corollary 2:}
Let $M^m$ be an entire spacelike graphic translator in $\mathbb{R}^{m+n}_n$ as defined in Theorem 5.  If there exists a number $\epsilon>0$ such that  the induced metric $(g_{ij})$ satisfies
\begin{equation}
\det(g_{ij})>\frac{\epsilon}{|x|} ,\;\;\;as\,\,\,\;\; |x|\rightarrow\infty , \end{equation}
then $M$ is an affine $m$-plane.

%%%%%%%%%%%%%%%%%%%%%%%%%%%%%%%%%%%%%%%%%%%%%%%%%%
\section{Preliminaries}
The pseudo-Euclidean space $\mathbb{R}^{m+n}_{n}$ is the linear space $\mathbb{R}^{m+n}$ endowed with the metric
\begin{equation}\label{2.1}
ds^{2} =\sum^{m}_{i=1}(dx^{i})^{2}-\sum^{m+n}_{\alpha=m+1}(dx^{\alpha})^{2}.
\end{equation}
Let $X:M^m\rightarrow\mathbb{R}^{m+n}_{n}$ be a spacelike $m$-submanifold in $\mathbb{R}^{m+n}_{n}$ with the second fundamental form $B$ defined by
\begin{equation}\label{2.2}
B_{YW}:=(\overline{\nabla}_{Y}W)^{\bot}
\end{equation}
for $Y,W\in\Gamma(TM)$, where $\overline \nabla$ denotes the connection on $\mathbb{R}^{m+n}_{n}$. Let $(\cdot)^{\top}$ and $(\cdot)^{\bot}$ denote the orthogonal projections into the tangent bundle $TM$ and the normal bundle $NM$, respectively. Let $\nabla$ and $\nabla^\bot$ be connections on the tangent bundle and the normal bundle of $M$, respectively.
Choose a local Lorentzian frame field $\{e_{i},e_{\alpha}\}(i=1,\cdots,m;\alpha=m+1,\cdots,m+n)$ such that $\{e_{i}\}$ are tangent vectors to $M$.
The mean curvature vector of $M$ in $\mathbb{R}^{m+n}_{n}$ is defined by
\begin{equation}\label{2.3}
H=\sum_\alpha H^{\alpha}e_{\alpha}=\sum_i B_{ii}=-\sum_{i,\alpha} h^{\alpha}_{ii}e_{\alpha},
\end{equation}
where $ h^{\alpha}_{ii}=\langle B_{ii}, e_\alpha\rangle$. Here, $\langle\cdot,\cdot\rangle$ is the canonical inner product in $\mathbb{R}^{m+n}_n$. We have the Gauss equation
\begin{equation}\label{2.4}R_{ijkl}=-\sum_\alpha(h^{\alpha}_{ik}h^{\alpha}_{jl}
-h^{\alpha}_{il}h^{\alpha}_{jk}),\end{equation}
and the Ricci curvature
\begin{equation}\label{2.5}R_{ij}=-\sum_{\alpha,k}(h^{\alpha}_{kk}h^{\alpha}_{ij}
-h^{\alpha}_{ki}h^{\alpha}_{kj}).\end{equation}

For a spacelike translating soliton $M^{m}$, by definition we can decompose the translating vector $T$ into a tangential part $V$ and a normal part $H$ on $M^{m}$, namely $T=V+H$.
Define $$\|H\|^{2}=-\langle H,H\rangle=-|H|^2,$$
where $\|H\|^{2}$ is absolute value of the norm square of the mean curvature. Similarly define $$\|B\|^{2}=-\langle B,B\rangle=-|B|^2.$$
From (\ref{2.3}) and the inequality of Schwartz, we have
\begin{equation}\label{2.6}
\|B\|^{2}\geq\frac{1}{m}\|H\|^{2}.
\end{equation}
Note that when the spacelike manifold $M^{m}$ is a Lagrangian gradient graph $(x,\nabla f)$ in $\mathbb{R}^{2n}_{n}$ with the indefinite metric $\sum_{i=1}^n dx_idy_i$, the functions $\|H\|^{2}$ and $\|B\|^{2}$ are the norm of Tchebychev vector field $\Phi$ and  the Pick invariant $J$ in relative geometry respectively (see\cite{X-H}), up to a constant. Therefore we can use some affine technique to estimate functions $\|H\|^{2}$ and $\|B\|^{2}$.

Let $z=\langle X,X\rangle$ be the pseudo-distance function on $M$. Then we have (see also \cite{X-Y-L-3})
\begin{equation}\label{2.7}
z_{,i} =e_{i}(z)=2\langle X,e_{i}\rangle,
\end{equation}
\begin{equation}\label{2.8}
z_{,ij}=Hess(z)(e_{i},e_{j})=2(\delta_{ij}-\langle X,h^{\alpha}_{ij}e_{\alpha}\rangle ),
\end{equation}
\begin{equation}\label{2.9}
\Delta z=\sum z_{,ii}=2m+2\langle X,H\rangle.
\end{equation}

In the following we set up the basic notations and formula for  an $m$-dimension spacelike graphic submanifold in $\mathbb{R}^{m+n}_n$. Let
$$M:=\{(x_{1},\cdots,x_{m},u^{1},\cdots,u^{n});x_{i}\in\mathbb{R},u^\alpha (x)=u^\alpha(x_1,\cdots,x_m)\},$$
where $i=1,\cdots,m$ and $\alpha=1,\cdots,n$. Denote
$$x=(x_{1}, \cdots,x_{m})\in\mathbb{R}^m;u=(u^{1} ,\cdots,u^{n})\in\mathbb{R}^n.$$
Let ${E_A} (A=1,\cdots,m+n)$ be the canonical Lorentzian basis of $\mathbb{R}^{m+n}$.
Namely, every component of the vector $E_A$ is 0, except that the $A$-th component is 1.
Then $$e_i = E_i +\sum_\alpha u^{\alpha}_i E_{m+\alpha},\quad  \,i \in\{1,\cdots,m\}$$
give a tangent frame on $M$. Here, $u^{\alpha}_i=\frac{\partial u^{\alpha}}{\partial x_i}$. In pseudo-Euclidean space with index $n$, the induced metric on spacelike submanifold $M$ is
$$g_{ij}=\langle e_i,e_j\rangle=\delta_{ij}-\sum_{\alpha}u^{\alpha}_iu^{\alpha}_j.$$
Then there are $n$ linear independent unit normal vectors,
$$e_{\alpha}=\frac{\sum_{i} u^{\alpha}_i E_i+E_{m+\alpha}}{(1-|Du^\alpha|^2)^\frac{1}{2}},\quad \,\alpha \in\{1,\cdots,n\},$$
where $Du^\alpha=(u^\alpha_1,\cdots,u^\alpha_m)$.
Thus the Levi-Civita connection with respect to the induced metric has the Chistoffel symbols
\begin{equation}\label{2.10}
\Gamma^{k}_{ij}=\frac{1}{2}g^{kl}(\frac{\partial g_{il}}{\partial x_{j}}+\frac{\partial g_{jl}}{\partial x_{i}}-\frac{\partial g_{ij}}{\partial x_{l}})=-g^{kl}u^{\alpha}_{ij}u^{\alpha}_l.
\end{equation}
Put a
nonzero constant vector
$$ T:=\sum a^{i}E_i+\sum b^{\alpha}E_{m+\alpha}.
$$
By the definition of the translator (\ref{1.2}), for each $e_\alpha$, there holds
$$\langle T,e_{\alpha}\rangle=-H^\alpha.$$
By calculation, we have
\begin{equation*}
\begin{aligned}
\langle H,e_\alpha\rangle&=g^{ij}\langle\overline{\nabla}_{e_i}{e_j},e_{\alpha}\rangle
=g^{ij}\langle u^{\beta}_{ij}E_{m+\beta},e_{\alpha}\rangle\\
&=\frac{-g^{ij}u^\alpha_{ij}}{(1-|Du^\alpha|^2)^\frac{1}{2}},
\end{aligned}
\end{equation*}
and
$$\langle T,e_{\alpha}\rangle=\frac{a^{i}u^{\alpha}_{i}-
 b^\alpha}{(1-|Du^\alpha|^2)^\frac{1}{2}}.$$
 Then (\ref{1.2}) is equivalent to the following elliptic system
\begin{equation}\label{2.12}
g^{ij}u^{\alpha}_{ij}=-a^{i}u ^{\alpha}_{i}+b^{\alpha}, \quad  \alpha \in\{1,\cdots,n\}.
\end{equation}

%%%%%%%%%%%%%%%%%%%%%%%%%%%%%%%%%%%%%%%%%%%%%%%
\section{Calculation of $\Delta \|H\|^{2}$ and $\Delta \|B\|^2$  }

Using the idea and technique in \cite{C-Q} and \cite{X-Y-L-3}, we have the following propositions.

\vskip 0.1in\textbf{Proposition 3.1:}
For the spacelike translating soliton $M^m$ in $\mathbb{R}^{m+n}_{n}$, the following estimate holds
\begin{equation*}
\Delta\|H\|^{2}\geq  \frac{2}{m} \|H\|^{4}- \langle T,\nabla\|H\|^{2}\rangle.
\end{equation*}

\vskip 0.1in\indent \textbf{Proof:}
Let $\{e_i\}$ be a local orthonormal normal frame field at the considered point of $M$. From (1.2), we derive
\begin{equation}\label{3.1}
\nabla_{e_j}^\bot H=(\overline{\nabla}_{e_j}H)^\bot=\left( \overline{\nabla}_{e_j}\left(T-\sum_k\langle T,e_k\rangle e_k\right)\right)^\bot=-\sum_k\langle T,e_k\rangle B_{jk},
\end{equation}
and
\begin{equation}\label{3.2}
\nabla_{e_i}^\bot\nabla_{e_j}^\bot H= -\sum_k\langle T,e_k\rangle\nabla_{e_i}^\bot B_{jk}-\sum_k\langle H,B_{ik}\rangle B_{jk}.
\end{equation}
Hence, using the Codazzi equation $$\nabla_{e_i}^\bot B_{jk}=\nabla_{e_k}^\bot B_{ji},$$ we have
\begin{equation}\label{3.3}
\begin{aligned} % requires amsmath; align* for no eq. number
\Delta|H|^{2}&=2|\nabla^\bot H|^{2}+2\langle H,\Delta ^\bot H\rangle \\
        & =2|\nabla^\bot H|^{2}-2 \langle H,\nabla_{V}^\bot H\rangle-2\sum_{i,k}\langle H,B_{ik}\rangle^{2}.
         \end{aligned}
\end{equation}
It follows that
\begin{equation}\label{3.4}
\begin{aligned}
\Delta\|H\|^{2}&=2\|\nabla^\bot H\|^{2}+2 \langle H,\nabla_{V}^\bot H\rangle+2\sum_{i,k}\langle H,B_{ik}\rangle^{2}\\
                              &\geq 2 \langle H,\nabla_{V} ^\bot H\rangle+\frac{2}{m} \|H\|^{4}.
 \end{aligned}
\end{equation}
Note that
\begin{equation}\label{3.5}
\begin{aligned}
2\langle H,\nabla_{V}^\bot H\rangle=\nabla_{V}\langle H,H\rangle=-\nabla_{V}\|H\|^{2}=-\langle T,\nabla\|H\|^{2}\rangle.
 \end{aligned}
\end{equation}
Then (\ref{3.4}) and (\ref{3.5}) together give proposition 3.1. \hfill$\Box$
\vskip 0.1in\noindent

In order to prove the completeness of entire graph spacelike translator with respect to the induced metric,
we need the following type estimate for $\Delta\|B\|^2$.

\vskip 0.1in\indent \textbf{Proposition 3.2:}
For the spacelike translating soliton $M$ in $\mathbb{R}^{m+n}_{n}$,
the Laplacian of the second fundamental form $B$ satisfies
\begin{equation*}
\Delta\|B\|^{2}\geq \frac{2}{n}\|B\|^{4}-\langle T,\nabla\|B\|^{2}\rangle.
\end{equation*}
\vskip 0.1in {\bf Proof:}
Let $\{e_i\}$ be a local tangent orthonormal frame field on $M$ and $\{e_\alpha\}$
a local normal orthonormal frame field on $M$ such that $\nabla e_{\alpha} = 0$ at the considered point.
From \cite{X-Y-L-3}, we have
\begin{equation}\label{3.6}
\begin{aligned}
\Delta\|B\|^{2}=\Delta\sum_{\alpha,i,j}(h^{\alpha}_{ij})^{2}=
2[\sum(h^{\alpha}_{ij,k})^{2}+h^{\alpha}_{ij}h^{\alpha}_{kk,ij}
-h^{\alpha}_{ij}h^{\alpha}_{li}h^{\beta}_{lj}h^{\beta}_{kk}\\ +h^{\alpha}_{ij}h^{\beta}_{ij}h^{\alpha}_{lk}h^{\beta}_{lk}
+(h^{\alpha}_{ij}h^{\beta}_{ik}-h^{\beta}_{lj}h^{\alpha}_{kl})
(h^{\alpha}_{qj}h^{\beta}_{qk}-h^{\beta}_{qj}h^{\alpha}_{kq})].
 \end{aligned}
\end{equation}
Let
$$S_{\alpha\beta}= \sum_{i,j}h^{\alpha}_{ij}h^{\beta}_{ij},~~~ N(h^{\alpha})=\sum_{i,j}(h^{\alpha}_{ij})^{2}, $$
then
$$\|B\|^{2}=\sum_{\alpha}S_{\alpha\alpha}.$$
So (\ref{3.6}) becomes
\begin{equation}\label{3.7}
\begin{aligned}
\Delta\|B\|^{2}=2[\sum(h^{\alpha}_{ij,k})^{2}- h^{\alpha}_{ij}H^{\alpha}_{,ij}+ h^{\alpha}_{ij}h^{\alpha}_{li}h^{\beta}_{lj}H^{\beta}\\ +\sum_{\alpha,\beta}(S_{\alpha\beta})^{2}
+N(h^{\alpha}h^{\beta}-h^{\beta}h^{\alpha})].
 \end{aligned}
\end{equation}
Note that
$$\sum_{\alpha,\beta}(S_{\alpha\beta})^{2}\geq
\frac{1}{n}(\sum_{\alpha}S_{\alpha\alpha})^{2}
=\frac{1}{n}\|B\|^{4},$$
$$N(h^{\alpha} h^{\beta}- h^{\beta} h^{\alpha})\geq0.$$
From (\ref{3.2}) and Codazzi equation, we get
\begin{equation*}
\begin{aligned}
\sum_{\alpha,i,j}h^{\alpha}_{ij}H^{\alpha}_{,ij}
=\sum_{i,j}\langle B_{ij},\nabla_{e_{i}}^{\bot}\nabla_{e_{j}}^{\bot} H\rangle
&=-\sum\langle B_{ij},\nabla_{V}^\bot B_{ij}\rangle-\sum_{i,j,k}\langle H,B_{ik}\rangle\langle B_{ij},B_{jk}\rangle\\
&=-\frac{1}{2}\langle V,\nabla\|B\|^2\rangle+\sum h^{\alpha}_{ij}h^{\alpha}_{jk}h^{\beta}_{ik}H^{\beta}.
\end{aligned}
\end{equation*}
Then substituting these inequalities into (\ref{3.7}) completes the proof of  proposition 3.2.\hfill$\Box$

%%%%%%%%%%%%%%%%%%%%%%%%%%%%%%%%%%%%%%%%%%

\section{Proof of theorem 3}

To gain Bernstein theorem, we are about to prove $\|H\|^{2}\equiv0 $ on $M$. If $\|H\|^{2}\neq 0$, then there exists a point $p_{0}\in M$ such that $\|H\|^{2}(p_{0})>0.$  Set $\|H\|^{2}(p_{0})=\lambda$. Denote by $r(p_{0},p)$  the geodesic distance function from $p_{0}\in M$  with respect to the induced metric $g$.  For any positive number $a$, let $B_{a}(p_{0}):=\{p\in M\;|\;r(p_{0},p)\leq a\}$. We consider the function
\begin{equation}\label{4.1}
\Psi:=(a^{2}-r^{2})^{2}\|H\|^{2}
\end{equation}
defined on $B_{a}(p_{0})$. Obviously, $\Psi$ attains its supremum at some interior point $p^{\ast}$. We may assume that $r^{2}$ is a $C ^{2}$-function in a neighborhood of $p^{\ast}$. Choose an orthonormal frame field on $M$ around $p^{\ast}$ with respect to the metric $g$. Then, at $p^{\ast}$, we have
$$\nabla\Psi=0,\quad\quad\Delta\Psi\leq 0.$$
Hence
\begin{equation}\label{4.2}
-\frac{2(r^{2})_{,i}}{a^{2}-r^{2}}+\frac{(\|H\|^{2})_{,i}}{\|H\|^{2}}=0,
\end{equation}
\begin{equation}\label{4.3}
-\frac{2\Delta r^{2}}{a^{2}-r^{2}}-\frac{2|\nabla r^{2}|^{2}}{(a^{2}-r^{2})^{2}}+\frac{\Delta\|H\|^{2}}{\|H\|^{2}}
-\frac{(\nabla\|H\|^{2})^{2}}{\|H\|^{4}}\leq0,
\end{equation}
where $","$  denotes the covariant derivative with respect to the metric $g$.\\
Inserting Proposition 3.1 into (\ref{4.3}), we get
\begin{equation}\label{4.4}
-\frac{2\Delta r^{2}}{a^{2}-r^{2}}-\frac{2|\nabla r^{2}|^{2}}
{(a^{2}-r^{2})^{2}}- \frac{\langle T,\nabla\|H\|^{2}\rangle}{\|H\|^{2}}
+\frac{2}{m}\|H\|^{2}-\frac{(\nabla\|H\|^{2})^{2}}{\|H\|^{4}}\leq0.
\end{equation}
Combining (\ref{4.2}) with (\ref{4.4}), we have
\begin{equation}\label{4.5}
0\geq-\frac{2\Delta r^{2}}{a^{2}-r^{2}}-\frac{24r^{2}}{(a^{2}-r^{2})^{2}}
 -\frac{\langle T,\nabla\|H\|^{2}\rangle}{\|H\|^{2}}+\frac{2}{m} \|H\|^{2}.
\end{equation}
By (\ref{4.2}) and the inequality of Schwarz, we obtain
\begin{equation}\label{4.6}
\begin{aligned}
 \langle T,\frac{\nabla\|H\|^{2}}{\|H\|^{2}}\rangle &=\langle V,\frac{\nabla\|H\|^{2}}{\|H\|^{2}}\rangle\\
&\leq\epsilon|V|^{2}+\frac{4 }{\epsilon}\frac{r^{2}}{(a^{2}-r^{2})^{2}}\\
&=\epsilon(C_0+\|H\|^{2})+\frac{4 }{\epsilon}\frac{r^{2}}{(a^{2}-r^{2})^{2}},
 \end{aligned}
\end{equation}
where $C_0=\langle T,T\rangle$ and $\epsilon$ is a small positive constant to be determined later. Inserting (\ref{4.6}) into (\ref{4.5}), we have
\begin{equation}\label{4.7}
\begin{aligned}
0&\geq%\frac{-2\Delta r^{2}}{a^{2}-r^{2}}-\frac{24r^{2}}{(a^{2}-r^{2})^{2}}+\frac{2\|\nabla H\|^{2}}{\|H\|^{2}}+2m\|H\|^{2}-\epsilon(c+\|H\|^{2})-
%\frac{4m^2}{\epsilon}\cdot\frac{r^{2}}{(a^{2}-r^{2})^{2}}\\
 -\frac{2\Delta r^{2}}{a^{2}-r^{2}}-(24+\frac{4 }{\epsilon})
 \cdot\frac{r^{2}}{(a^{2}-r^{2})^{2}}+(\frac{2}{m} -\epsilon)\|H\|^{2}-\epsilon C_0.
\end{aligned}
\end{equation}

Now we shall use the technique in \cite{J-L} to calculate the term $\Delta r^{2}$. In the case $p_{0}\neq p^{\ast}$, we denote $a^{\ast}=r(p_{0},p^{\ast})>0$ and assume that
$$\max_{B_{a^{\ast}}(p_{0})}\|H\|^{2}=\|H\|^{2}(\tilde{p}).$$
By Proposition 3.1 we have
$$\max_{B_{a^{\ast}}(p_{0})}\|H\|^{2} =\max_{\partial B_{a^{\ast}}(p_{0})}\|H\|^{2}.$$
Thus
\begin{equation}\label{4.8}
   (a^{2}-r^{2}(p^*))^{2}\|H\|^{2}(\tilde p)  =   (a^{2}-r^{2}(\tilde p))^{2}\|H\|^{2}(\tilde p)\leq (a^{2}-r^{2}(p^*))^{2}\|H\|^{2}(p^*).
\end{equation}
For any $p\in B_{a^{\ast}}(p_{0})$, by the Gauss equation we know that the Ricci curvature $Ric(M,g)$ is bounded from below by
\begin{equation}\label{4.9}
\begin{aligned}
R_{ii}(p)&\geq-\frac{1}{4} \|H\|^{2}(\tilde{p}).
 \end{aligned}
\end{equation}
By the Laplacian comparison theorem, we get
\begin{equation}\label{4.10}
\begin{aligned}
\Delta r^{2}%&\leq2[1+(m-1)(1+\sqrt{-K}r)]\\
            %&=2[1+(m-1)(1+\frac{1}{2}m\|H\|r)]\\
            &\leq2m+(m-1) r\|H\|(\tilde{p}).
 \end{aligned}
\end{equation}
Substituting (\ref{4.10}) into (\ref{4.7}), we have
\begin{equation}\label{4.11}
 (\frac{2}{m} -\epsilon)\|H\|^{2}(p^{\ast})
 \leq c(m,\epsilon)\frac{a^{2}}{(a^{2}-r^{2})^{2}}
+ \epsilon C_0+2(m-1) \|H\|(\tilde{p})\frac{a}{a^{2}-r^{2}},
 \end{equation}
where $c(m,\epsilon)$ is a positive constant depending only on $m$ and $\epsilon$.
Multiplying both sides of (4.11) by $(a^2-r^2)^2 (p^*)$, we have
\begin{equation}\label{4.12}
(\frac{2}{m} -\epsilon)\Psi
 \leq c(m,\epsilon) a^{2}
+\epsilon C_{1} a^4 +2(m-1)  a(a^{2}-r^{2})(p^*)\|H\|(\tilde{p}),
\end{equation}
where
$$ C_{1} = \left\{ \begin{array}{ll}
C_0, & \mbox{if}\  C_0\geq0 \\
0, & \mbox{if}\ C_0<0.
\end{array}
\right.
$$
By (4.8) and the inequality of Schwarz, we get
\begin{equation}\label{4.13}
 \Psi\leq c(m,\epsilon) a^{2}
+\frac{mC_{1}\epsilon}{2 -2m\epsilon} a^4.
\end{equation}
In the case $p_{0}= p^{\ast}$, we have $r(p_{0},p^{\ast})=0$. It is easy to get (4.13) from (4.7). Therefore (4.13) holds on $B_{a}(p_{0})$.

Then, at any interior point $q\in B_{a}(p_{0})$, we get
\begin{equation}\label{4.14}
\begin{aligned}
\|H\|^{2}(q)&\leq  c(m,\epsilon) \frac{a^{2}}{(a^{2}-r^{2})^{2}}+\frac{mC_{1}\epsilon}{2 -2m\epsilon}\frac{a^{4}}{(a^{2}-r^{2})^{2}}.
 \end{aligned}
\end{equation}
\textbf{Case 1.} If $C_{0}>0$, we choose $\epsilon$ sufficiently small, such that
$$\frac{mC_{1}\epsilon}{2 -2m\epsilon}\leq\frac{\lambda}{2}.$$
For $a\rightarrow\infty$, at the point $q$, we get
$$\|H\|^{2}\leq\frac{2\lambda}{3}.$$
In particular,
$$\|H\|^{2}(p_{0})\leq\frac{2\lambda}{3}.$$
This contradicts to $\|H\|^{2}=\lambda$ at $p_{0}$, so we can conclude that $\|H\|^{2}\equiv0$ on $M^n$. Therefore it is an affine $m$-plane by theorem 1.
\\
\textbf{Case 2.} If $C_{0}= 0$, let $a\rightarrow\infty$ in (\ref{4.14}), we get
$$\|H\|^{2}=0.$$
On the other hand, it is easy to see that $H$ is a zero vector, since $M$ is an $m$-dimensional spacelike submanifold. Thus $T\in TM$. This contradicts to the translating vector $T$ is nonzero. Then there exists no complete $m$-dimensional spacelike translating soliton  with a \textbf{lightlike} translating vector.
\\
\textbf{Case 3.} If $C_{0}<0$, let $a\rightarrow\infty$ in (\ref{4.14}), we get
$$\|H\|^{2}\leq0.$$
This is impossible. Then there exists no complete $m$-dimensional spacelike translating soliton  with a \textbf{timelike} translating vector.
This completes the proof of Theorem 3.
\hfill $\Box$

\vskip 0.1in\textbf{Remark}: Here we mention that we can replace $\|H\|$
 with $\|B\|$ in the proof of theorem 3 and use proposition 3.2 to prove the second fundamental form $\|B\|\equiv 0$. It directly prove $M^m$ must be a plane.

\section{Proof of theorem 5}

To gain theorem 5, we will show that the graph spacelike translating soliton is complete with respect to the induced metric $g$. Thus from theorem 3, we complete the proof of theorem 5. Using the similar calculation of section 5 in \cite{X-Z}, we prove the completeness of the induced metric $g$ if $\|B\|$ has an upper bound. To the simplicity, here we shall use theorem 2.1 of Prof. Xin  in \cite{X-Y-L-2} to obtain the completeness.

\vskip 0.1in  \textbf{ Proposition 5.1:}
Let $u^\alpha$  be smooth functions defined everywhere in $\mathbb{R}^m$. Suppose their graph $M=(x,u(x))$ is a spacelike translator in $\mathbb{R}^{m+n}_n$. If the norm of the second fundamental form $\|B\|$ has a bound, then $M$ is complete with respect to the induced metric.

\vskip 0.1in\indent \textbf{Proof:}
 Without loss of generality, we assume that the origin $0\in M$. From Proposition 3.1 of \cite{J-X}, we know that the pseudo-distance function $z=\langle X,X\rangle$ on $M$ is a non-negative proper function. By (\ref{2.6}), we know that if $\|B\|$ has an upper bound, $\|H\|$ also has an upper bound. On the other hand, by (3.1) we get
 \begin{equation*}
\begin{aligned}
\|\nabla^\bot H\|^2&=-\langle \nabla^\bot H,\nabla^\bot H \rangle=-\sum_j\langle\sum_k\langle T,e_k\rangle B_{jk},\sum_l\langle T,e_l\rangle B_{jl} \rangle\\
&\leq|V|^2\|B\|^2\leq (|T|^2+\|H\|^2)\|B\|^2.
 \end{aligned}
\end{equation*}
Then $\|\nabla^\bot H\|$ also has a bound.  By Theorem 2.1 in \cite{X-Y-L-2}, we have
for some $k>0$, the set $\{z\leq k\}$ is compact, then there is a constant $c$ depending only on
the dimension $m$ and the bounds of mean curvature and its covariant derivatives
such that for all $x\in M$ with $z(x)\leq \frac{k}{2}$,
\begin{equation}\label{5.1}
|\nabla z|\leq c(z+1).
\end{equation}
Let $\gamma:[0,r]\rightarrow M$ be a geodesic on $M$ issuing from the origin 0. Integrating
(\ref{5.1}) gives
$$z(\gamma (r))+1\leq \exp(cr),$$
which forces $M$ to be complete with respect to the induced metric.
 \hfill $\Box$

\vskip 0.1in In the following we prove the bound of $\|B\|$.

\vskip 0.1in  {\bf Proposition 5.2:}
  Under the assumption in Theorem 5, the function $\|B\|$ has an upper bound on $\mathbb{R}^{m}$.

\vskip 0.1in  \textbf{Proof:}
By assumption in theorem 5, there exists a large enough $R_0$ such  that
$(g_{ij})>\frac{\epsilon}{|x|}I$ in case $|x|>R_0$. Choose a positive number $a>R_0$.
Let $B_a(0):=\{x\in \mathbb{R}^m\;|\; |x|\leq a\}$. Consider the  function below
\begin{equation*}
F(x):=  (a^2-|x|^2)^{2} \|B\|^2
 \end{equation*}
defined on $B_a( 0 )$. Obviously,  $F$ attains its supremum at
some interior point $p^*$. We can assume that $ \|B\|(p^*)>0$.  Then, at $p^*$,
$$ (\log F)_i=0, \quad (\log F)_{ij}\leq 0.$$
In the following we calculate the first equation explicitly,
and contract the second term above with the positive definite matrix $(g^{ij})$.
\begin{equation}\label{5.2}
\frac{(\|B\|^2)_{i}}{\|B\|^2} - \frac{4x_{i}}{a^2-|x|^2}=0,
\end{equation}
\begin{equation}\label{5.3}
\sum_{i,j=1}^{n} g^{ij}\left( \frac{
(\|B\|^2)_{ij}}{\|B\|^2}-\frac{ (\|B\|^2)_{i} (\|B\|^2)_{j}}{\|B\|^4} - \frac{ 8 x_ix_j
}{(a^2-|x|^2)^2} -\frac{4 \delta_{ij}}{a^2-|x|^2}\right)   \leq
0.\end{equation}
From Proposition 3.2 and the inequality of Schwarz, we get
\begin{equation}\label{5.4}
\begin{aligned}
\sum_{i,j} g^{ij}\frac{(\|B\|^2)_{ij}}{\|B\|^2}&=\frac{\Delta\|B\|^2}{\|B\|^2}+\sum g^{ij}\Gamma^{k}_{ij}\frac{(\|B\|^2)_{k}}{\|B\|^2}\\
&\geq\frac{2}{n}\|B\|^2- \langle V,\frac{\nabla\|B\|^2}{\|B\|^2}\rangle+\sum_{i,j,k}g^{ij}
\Gamma^{k}_{ij}\frac{(\|B\|^2)_{k}}{\|B\|^2}\\
&\geq\frac{2}{n}\|B\|^2- \frac{1}{mn}|V|^2-mn\frac{(\nabla\|B\|^2)^2}{\|B\|^4}
+\sum_{i,j,k}g^{ij}\Gamma^{k}_{ij}\frac{(\|B\|^2)_{k}}{\|B\|^2}.
 \end{aligned}
\end{equation}
%%%%%%%%%%%%%%%%%%%%%%%%%%%%%%%%%%

In the following we shall estimate the term $\sum_{i,j,k}g^{ij}\Gamma^{k}_{ij}\frac{(\|B\|^2)_{k}}{\|B\|^2}$. By (\ref{2.10})  (\ref{2.12}) and (\ref{5.2}), we have
\begin{equation*}
\sum_{i,j,k}g^{ij}\Gamma^{k}_{ij}\frac{(\|B\|^2)_{k}}{\|B\|^2}
 =\sum_{\alpha}(a^{i}u^{\alpha}_{i}-b^{\alpha})
g^{kl}u^{\alpha}_{l}\frac{4x_k}{a^2-|x|^2}.
 \end{equation*}
By the inequality of Cauchy,  we get
\begin{equation*}
\sum_{\alpha}(a^{i}u^{\alpha}_{i}-b^{\alpha})
g^{kl}u^{\alpha}_l x_k
\leq \left[\sum_{\alpha}(a^{i}u^{\alpha}_{i}-b^{\alpha})^2\right]^{\frac{1}{2}}
\left[\sum_{\alpha}(g^{kl}u^{\alpha}_{l}x_{k})^2\right]^{\frac{1}{2}}.
 \end{equation*}
By the spacelike of \ $M$,  we have$$\sum_\alpha(u^\alpha_i)^2<1,\quad\quad \forall\quad  1\leq i\leq m.$$
Then
 \begin{equation*}
 \begin{aligned}
\left[\sum_{\alpha}(a^{i}u^{\alpha}_{i}-b^{\alpha})^2\right]^{\frac{1}{2}}
\leq &  \left[\sum_{\alpha}(b^{\alpha})^2\right]^{\frac{1}{2}}
+\left [\sum_{\alpha}(\sum_{i}a^{i}u^{\alpha}_{i})^2\right ]^{\frac{1}{2}}\\
\leq & |\vec{b}|+ \left[|\vec{a}|^2(\sum_{\alpha}\sum_{i}(u^{\alpha}_{i})^2)\right]^{\frac{1}{2}}\\
 = &|\vec{b}|+|\vec{a}|\left[\sum_{i}(\sum_{\alpha}u^{\alpha}_{i})^2
 \right]^{\frac{1}{2}}\\
 \leq &|\vec{b}|+\sqrt{m}|\vec{a}|,
  \end{aligned}
 \end{equation*}
where \ $\vec{a}=(a^1,...,a^m)$ ,\ $\vec{b}=(b^1,...,b^n)$. On the other hand, \\
 \begin{equation*}
 \begin{aligned}
 \left[\sum_{\alpha}(g^{kl}u^{\alpha}_{l}x_{k})^2\right]^{\frac{1}{2}}
 \leq &\left[\sum_{\alpha}(g^{kl}u^{\alpha}_{l}u^{\alpha}_{k})\right]
 ^{\frac{1}{2}}(g^{kl}x_{l}x_{k})^{\frac{1}{2}}\\
 \leq &\left[(\sum_{i}g^{ii})(\sum_{\alpha}
(\sum_{k}(u^{\alpha}_{k})^2))\right]^{\frac{1}{2}}
   \left[(\sum_{k}g^{kk})| {x}|^2\right]^{\frac{1}{2}}\\
    \leq &\left[(\sum_{i}g^{ii})(\sum_{k}
(\sum_{\alpha}(u^{\alpha}_{k})^2))\right]^{\frac{1}{2}}
    (\sum_{k}g^{kk}) ^{\frac{1}{2}}| {x}|\\
    \leq &\sqrt{m}  a (\sum_{i}g^{ii}).
     \end{aligned}
 \end{equation*}
Therefore
\begin{equation}\label{5.5}
\sum_{i,j,k}g^{ij}\Gamma^{k}_{ij}\frac{(\|B\|^2)_{k}}{\|B\|^2}
\leq\frac{C_1a\sum g^{ii}}{a^2-|x|^2},
\end{equation}
where \ $C_1=4m (|\vec{a}|+|\vec{b}|)$.

Thus, inserting (\ref{2.6}), (\ref{5.2}) and (\ref{5.5}) into (\ref{5.4}), we get
\begin{equation}\label{5.6}
\begin{aligned}
\sum_{i,j} g^{ij}\frac{(\|B\|^2)_{ij}}{\|B\|^2}
&\geq\frac{2}{n}\|B\|^2- \frac{1}{mn}|V|^2-mn\frac{(\nabla\|B\|^2)^2}{\|B\|^4}-\frac{C_1a\sum g^{ii}}{a^2-|x|^2}\\
& \geq \frac{1}{n} \|B\|^2-\frac{1}{mn} C_0- 16mn\frac{ a^2\sum g^{ii}}{(a^2-|x|^2)^2}-\frac{C_1a\sum g^{ii}}{a^2-|x|^2},
  \end{aligned}
 \end{equation}
where $C_0=\langle T,T\rangle$.\\
Then, inserting (\ref{5.2}) and (\ref{5.6}) into (\ref{5.3}), we have
\begin{equation}\label{5.7}
 \|B\|^2 \leq \frac{1}{m} C_0+ (16mn^2+28n) \frac{a^2\sum g^{ii}}{(a^2-|x|^2)^2}+n C_1\frac{a\sum g^{ii}}{a^2-|x|^2}.
 \end{equation}
Multiplying both sides of (\ref{5.7}) by $(a^2-|x|^2)^2(p^*)$, we have
\begin{equation}\label{5.8}
\begin{aligned}
(a^2-|x|^2)^{2} \|B\|^2(p^*)
\leq \frac{1}{m}C_0a^{4}+[(16mn^2+28n)a^2+n C_1a^3]\sum g^{ii}.
\end{aligned}
 \end{equation}
If $p^*\in B_{R_0}(0)$, then for any $x\in B_{a}(0)$,
$$F(x)\leq a^4\max_{B_{R_0}(0)} \|B\|^2.$$
If $p^*\in B_a(0)\backslash B_{R_0}(0)$, by assumption  $$\sum g^{ii}\leq m\frac{|x|}{\epsilon},$$
we get
\begin{equation}\label{5.9}
 [(a^2-|x|^2)^{2} \|B\|^2] (p^*)
\leq \frac{1}{m} C_0a^{4}+ C_2 a^3+\frac{mn C_1}{\epsilon}a^4.
\end{equation}
where $C_2$ is a
positive constant depending only on $m$, $n$ and $\epsilon$.
Thus, in both cases, we always have
\begin{equation}\label{5.10}
F(x)\leq[(a^2-|x|^2)^{2} \|B\|^2] (p^*)
\leq C_3 a^4+C_2a^3,
\end{equation}
where $C_3$ is a
positive constant depending only on $m$, $n$, $\epsilon$, $C_1$, $C_0$  and  $\max_{B_{R_0}(0)}\|B\|^2$.
Then, at any interior point $q\in B_a(0)$, we obtain
\begin{equation}\label{5.11}
[(a^2-|x|^2)^2\|B\|^2](q)\leq [(a^2-|x|^2)^2\|B\|^2](p^*)\leq C_3 a^4+C_2a^3.
\end{equation}
Dividing both sides of (\ref{5.11}) by $(a^2-|x|^2)^2$, we obtain
\begin{equation}\label{5.12}
\|B\|^2(q) \leq C_3 \frac{a^4}{(a^2-|x|^2)^{2}}+C_2\frac{a^3}{(a^2-|x|^2)^{2}}.
\end{equation}
Let $a\rightarrow\infty$ in (\ref{5.12}), we get
$$\|B\|^2(q)\leq C_3.
$$ This completes the proof of Proposition 5.2.\hfill$\Box$

\vskip 0.05in\noindent\textbf{Acknowledgements}: We wish to express our sincere gratitude to Prof. Qun Chen and Xingxiao Li for thier valuable and helpful discussions on the topic.

%%%%%%%%%%%%%%%%%%%%%%%%%%%%%%%%%%
%%%%%%%%%%%%%%%%%%%%%%%%%%%%%%%%%%
%%%%%%%%%%%%%%%%%%%%%%%%%%%%%%%%%%%%%%%%%%%%%%
%%%%%%%%%%%%%%%%%%%%%%%%%%%%%%%%
%%%%%%%%%%%%%%%%%%%%%%%%%%%%%%%%%

\newpage

\noindent Ruiwei Xu\quad\quad Tao Liu\\
School of Mathematics and Information Science,\\
Henan Normal University,\\
Xinxiang, Henan 453007, P.R. China,\\
E-mail: rwxu$@$henannu.edu.cn \\
E-mail: liutao-319$@$163.com

%%%%%%%%%%%%%%%%%%%%%%%%%%
%%%%%%%%%%%%%%%%%%%%%%%%%%
%%%%%%%%%%%%%%%%%%%%%%%%%%
%%%%%%%%%%%%%%%%%%%%%%%%%%
%%%%%%%%%%%%%%%%%%%%%%%%%%
%%%%%%%%%%%%%%%%%%%%%%%%%%
\end{document}